

\magnification=1100
\overfullrule0pt

\input prepictex
\input pictex
\input postpictex
\input amssym.def


\def\qed{\hbox{\hskip 1pt\vrule width4pt height 6pt depth1.5pt \hskip 1pt}}
\def\mapright#1{\smash{
   \mathop{\longrightarrow}\limits^{#1}}}

\def\mapdown#1{\Big\downarrow
   \rlap{$\vcenter{\hbox{$\scriptstyle#1$}}$}}

\def\CC{{\Bbb C}}
\def\FF{{\Bbb F}}
\def\GG{{\Bbb G}}
\def\PP{{\Bbb P}}
\def\QQ{{\Bbb Q}}
\def\RR{{\Bbb R}}
\def\ZZ{{\Bbb Z}}

\def\cB{{\cal B}}

\def\cE{{\cal E}}
\def\cF{{\cal F}}
\def\cI{{\cal I}}

\def\cL{{\cal L}}
\def\cO{{\cal O}}
\def\cT{{\cal T}}
\def\gb{{\goth b}}
\def\gg{{\goth g}}
\def\gh{{\goth h}}
\def\gp{{\goth p}}

\def\Im{{\rm Im}}
\def\ad{{\rm ad}}
\def\ch{{\rm ch}}
\def\gr{{\rm gr}}
\def\td{{\rm td}}




\font\smallcaps=cmcsc10
\font\titlefont=cmr10 scaled \magstep1

\font\sectionfont=cmbx10
\font\tinyrm=cmr10 at 8pt


\newcount\sectno
\newcount\subsectno
\newcount\resultno

\def\section #1. #2\par{
\sectno=#1
\resultno=0
\bigskip\noindent{\sectionfont #1.  #2}~\medbreak}

\def\subsection #1\par{\bigskip\noindent{\it  #1} \medbreak}


\def\prop{ \global\advance\resultno by 1
\bigskip\noindent{\bf Proposition \the\sectno.\the\resultno. }\sl}
\def\lemma{ \global\advance\resultno by 1
\bigskip\noindent{\bf Lemma \the\sectno.\the\resultno. }
\sl}
\def\remark{ \global\advance\resultno by 1
\bigskip\noindent{\bf Remark \the\sectno.\the\resultno. }}
\def\example{ \global\advance\resultno by 1
\bigskip\noindent{\bf Example \the\sectno.\the\resultno. }}
\def\cor{ \global\advance\resultno by 1
\bigskip\noindent{\bf Corollary \the\sectno.\the\resultno. }\sl}
\def\thm{ \global\advance\resultno by 1
\bigskip\noindent{\bf Theorem \the\sectno.\the\resultno. }\sl}
\def\defn{ \global\advance\resultno by 1
\bigskip\noindent{\it Definition \the\sectno.\the\resultno. }\slrm}

\def\endprop{\rm\bigskip}

\def\pf{\rm\bigskip\noindent{\it Proof. }}
\def\endpf{\qed\hfil\bigskip}


\def\formula{\global\advance\resultno by 1
\eqno{(\the\sectno.\the\resultno)}}
\def\formulano{\global\advance\resultno by 1 (\the\sectno.\the\resultno)}
\def\tableno{\global\advance\resultno by 1
\the\sectno.\the\resultno. }
\def\lformula{\global\advance\resultno by 1
\leqno(\the\sectno.\the\resultno)}

\def\monthname {\ifcase\month\or January\or February\or March\or April\or
May\or June\or
July\or August\or September\or October\or November\or December\fi}

\newcount\mins  \newcount\hours  \hours=\time \mins=\time
\def\now{\divide\hours by60 \multiply\hours by60 \advance\mins by-\hours
     \divide\hours by60         
     \ifnum\hours>12 \advance\hours by-12
       \number\hours:\ifnum\mins<10 0\fi\number\mins\ P.M.\else
       \number\hours:\ifnum\mins<10 0\fi\number\mins\ A.M.\fi}


\nopagenumbers
\def\runningtitle{\smallcaps a pieri-chevalley formula for k(g/b)}
\headline={\ifnum\pageno>1\eoheadline\else\firstheadline\fi}
\def\names{\smallcaps h. pittie\quad and\quad a. ram}
\def\firstheadline{\noindent Preliminary Draft \hfill  \today}
\def\firstheadline{}
\def\eoheadline{\ifodd\pageno\oddheadline\else\evenheadline\fi}
\def\oddheadline{\tenrm\hfil\runningtitle\hfil\folio}
\def\evenheadline{\tenrm \folio\hfil{\names}\hfil}

\vphantom{$ $}  
\vskip.75truein

\centerline{\titlefont A Pieri-Chevalley formula for K(G/B)}
\bigskip
\centerline{\rm H. Pittie}
\centerline{Department of Mathematics, Graduate Center}
\centerline{City University of New York}
\centerline{New York, NY 10036}
\bigskip
\centerline{\rm Arun Ram${}^\ast$}
\centerline{Department of Mathematics}
\centerline{Princeton University}
\centerline{Princeton, NJ 08544}
\centerline{{\tt rama@math.princeton.edu}}
\centerline{Preprint: May 19, 1998}

\footnote{}{\tinyrm ${}^\ast$ Research supported in part by National
Science Foundation grant DMS-9622985.}

\bigskip


\section 0. Introduction

Algebraic combinatorics and geometry come together in a beautiful
way in the study of the cohomology $H^*(G/B)$ of the generalized flag variety
$G/B$ for a complex semisimple Lie group $G$ with Borel subgroup $B$. 
The cohomology $H^*(G/B)$ is isomorphic (as a graded ring) to the quotient of a
polynomial ring by the ideal generated by $W$-symmetric functions without constant
term and has a natural basis of ``Schubert polynomials'' $[X_w]\in H^*(G/B)$
where $[X_w]$ is the Poincar\'e dual of the fundamental class of the Schubert
variety
$X_w\subseteq G/B$ in $H_*(G/B)$.  One of the fundamental results
in the theory of Schubert polynomials and the cohomology of the flag
variety is a formula of Chevalley which gives an expansion of the product
$\lambda\cdot [X_w]$ in terms of the Schubert class basis for an element
$\lambda\in H^2(G/B)$. 

A similar picture holds for the K-theory of $G/B$.  The ring $K(G/B)$ is
isomorphic to a quotient of a {\it Laurent} polynomial ring by an ideal generated
by certain $W$-symmetric functions and has a basis given by classes $[\cO_{X_w}]$
where $\cO_{X_w}$ is the structure sheaf of the Schubert variety $X_w\subseteq
G/B$ extended by $0$ outside of $X_w$.  In this paper we give an analogue
of Chevalley's formula for the ring $K(G/B)$.  Specifically, we give an
explicit combinatorial formula for $e^\lambda[\cO_{X_w}]$, the tensor product of a
(negative) line bundle with the structure sheaf of a Schubert variety, expanded in
terms of the Schubert class basis $\{[\cO_{X_w}]\}$.

The Chern character is an isomorphism between $K(G/B)$ and $H^*(G/B)$
and Chevalley's formula can be recovered from ours by 
applying the Chern character and comparing lowest degree terms.  
The higher order terms of our formula may yield
further interesting identities in cohomology.

Fulton and Lascoux [FL] have given a formula similar to ours for
the $G=GL_n(\CC)$ case.  Our formula is a generalization of their formula to
general type except that we work only with $K(G/B)$ instead of the K-theory of
the flag bundle.  In our work the column strict tableaux used by Fulton and
Lascoux are replaced by Littelmann's path model. This has two advantages: (1) it
allows us to work in general type and (2) it obviates the need for the complex
combinatorics associated with the jeu de taquin and the ``rectification'' of
tableaux.  It is possible that the more general flag bundle result of Fulton
and Lascoux can be obtained from our general commutation formula, Theorem 4.2, but
to sort this out properly one would have to understand concretely the connection
between the tableaux and the Littelmann paths.

In the first section of the paper we review notations and recall why one
is able to work with $K(G/B)$ as a quotient of a Laurent polynomial
ring.  In the second section we derive an expression for
$[\cO_p]\in K(G/P)$, for a point $p\in G/P$, in terms of familiar vector
bundles.  We are able to pull back this formula into $K(G/B)$ to obtain
formulas for certain special Schubert classes in terms of line bundles.
In section 3 we recall the operators which play the same role as the BGG
operators in cohomology and show that they can be used to give explicit
(inductive) expressions for the Schubert classes in $K(G/B)$.  In section 4 we
prove the main theorem, which gives a commutation relation between line bundles
and the Schubert classes. Our new Pieri-Chevalley formula is an immediate
consequence of this relation. In the final section we explain how the K-theory
relates to cohomology and how our formula implies the classical Chevalley formula.

The main results of the preliminary sections can all be considered well
known. These results can be found, either explicitly or implicitly, in the
work of Demazure [D], Kostant and Kumar [KK], Fulton and Lascoux [FL]
and others.  For the convenience of the reader, we have given short proofs or
sketches of proofs for most of these results and, in the final section,
we have given a dictionary between K-theory and $H^*(G/B)$.  This
dictionary illustrates how our results relate to the theory
of Schubert polynomials.

We would like to thank M. Green, S. Kumar, T. Shifrin, and A.
Vasquez for stimulating conversations during our work on this paper.

\section 1. Background

In this section we set up notations and recall how one is able
to work with the K-theory of $G/B$ ``purely combinatorially".
This is possible because (in analogy with the cohomology
$H^*(G/B)$ and the combinatorial theory of Schubert polynomials)
the K-theory of $G/B$ is isomorphic to a quotient of
$R(T)$, which is a ring of Laurent polynomials. 
In section 5 we shall explain how our results about $K(G/B)$ relate to
and generalize well known results about Schubert polynomials and $H^*(G/B)$.

\medskip
\subsection{The K-theory $K(G/B)$}

If $X$ is a quasi-projective
variety let
$$\eqalign{
K(X)&=\hbox{the Grothendieck group of coherent sheaves on $X$}, \cr
K_{lf}(X)&=\hbox{the Grothendieck group of locally free sheaves on $X$}, \cr
K_{vb}(X)&=\hbox{the Grothendieck group of vector bundles on $X$}. \cr
}$$
If $f\colon X\to Y$ is a morphism of projective varieties then we have maps
$$\matrix{
f^!\colon &K(Y) &\longrightarrow &K(X) \cr
&[\cF] &\longmapsto &[f^*\cF], \cr
\cr
f_!\colon &K(X) &\longrightarrow &K(Y) \cr
&[\cF] &\longmapsto &\displaystyle{ \sum_i (-1)^i[R^if_*\cF]. }\cr
}
$$
See [CG, \S 5.2] for K-theory background. 
If $X$ is smooth then the isomorphism between
$K(X)$  and $K_{lf}(X)$ is given by assigning to the class of a sheaf the
alternating sum of the sheaves in a locally free resolution, see [BS, \S4] and
[Ha, III Ex. 6.9].  There is always a map $K_{lf}(X)\to K_{vb}(X)$ which
assigns to a locally free sheaf the underlying vector bundle and the results of
[P] imply that this map is an isomorphism when $X=G/B$, see Proposition 1.5 below.
Thus, in the case which we wish to consider in this paper, $X=G/B$, all three
K-theories are isomorphic. 

Let $G$ be a complex connected simply connected semisimple Lie group.  Fix 
a maximal torus $T$ and a Borel subgroup $B$ such that 
$T\subseteq B\subseteq G$.  
The Bruhat decomposition says that $G$ is a disjoint union
of double cosets of $B$ indexed by the elements of the Weyl group $W$, 
$G = \bigcup_{w\in W} BwB$.
The {\it flag variety} is the projective variety formed by the coset space 
$G/B$ and the Bruhat decomposition of $G$ induces a cell decomposition of
$G/B$.  For each $w\in W$ the subset $X_w^\circ=BwB/B$ is the 
{\it Schubert cell} and its closure $X_w$ is the {\it Schubert variety}.
The formulas
$$\dim(X_w^\circ)=\ell(w)
\qquad\hbox{and}\qquad X_w = \bigcup_{v\le w} X_v^\circ \formula$$
define the {\it length} $\ell(w)$ of $w\in W$ and 
the {\it Bruhat-Chevalley order} $\le$ on the Weyl group, respectively.
It follows from the Bruhat decomposition  (see lecture 4 by Grothendieck in [C])
that 
$$\hbox{$K(G/B)$ is a free $\ZZ$-module with basis $\{ [\cO_{X_w}]\ |\ w\in W\}$,
}$$
where $X_w$ are the Schubert varieties in $G/B$ and $\cO_{X_w}$ is the structure
sheaf of $X_w$ extended to $G/B$ by defining it to be $0$ outside $X_w$.

\medskip
\subsection{The isomorphism $K(G/B)\cong R(T)/\cI$}

For any group $H$ let $R(H)$ be the Grothendieck group of complex representations
of $H$. Let $\Lambda = \sum_{i=1}^n \ZZ\omega_i,$ where the $\omega_i$ are 
the fundamental weights of the Lie algebra $\gg$ of $G$. 
We shall use the ``geometric'' convention (see [CG, 6.1.9(ii)]) and let
$e^{-\lambda}$ be the element of $R(T)$ corresponding to the character determined
by $\lambda\in \Lambda$.
Then 
$$\hbox{$R(T)$ has $\ZZ$-basis $\{e^\lambda\ |\ \lambda\in \Lambda\}$,
\qquad 
with multiplication 
$e^\lambda e^\mu = e^{\lambda+\mu}$,}$$
and Weyl group action determined by
$we^\lambda=e^{w\lambda}$, for $w\in W$ and $\lambda\in \Lambda$.
In this way $R(T)$ is a Laurent polynomial ring and $R(G)\cong R(T)^W$ is the
subalgebra of ``symmetric functions'' in $R(T)$. 

Suppose that $V$ is a $T$-module.  Since $T\cong B/U$, where $U$ is
the unipotent radical of $B$, we can extend $V$ to be a $B$-module by 
defining the action of $U$ to be trivial.  Define a vector bundle
$$\matrix{\pi\colon &G\times_B V &\longrightarrow &G/B \cr
&(g,v) &\longmapsto &gB \cr}
\qquad\hbox{where}\qquad
G\times_B  V= {G\times V\over \langle (g,v)\sim (gb,b^{-1}v) \rangle},$$
so that  $G\times_B V$ is the set of pairs $(g,v)$, $g\in G$, $v\in V$,
modulo the equivalence relation $(g,v)\sim (gb,b^{-1}v)$.
This construction induces a ring homomorphism
$$\matrix{
\phi\colon &R(T) &\longrightarrow &K(G/B) \cr
&V &\mapsto & (G\times_B V \mapright{\pi} G/B). \cr
}\formula$$
If $V$ is a $G$-module then $\phi(V)=\dim(V)$ in $K(G/B)$.
This is because the map
$$\matrix{
G\times_B V &\longrightarrow &G/B \times V\cr
(g,x) &\longmapsto &(gB,gx) \cr
}\formula$$
is an isomorphism between $G\times_B V$ and the trivial bundle
$G/B\times V$.
Define $\varepsilon \colon R(T)\to \ZZ$ by $\varepsilon(e^\lambda)=1$
for $\lambda\in \Lambda$.  Then the map $\phi$ in (1.2) gives an
isomorphism (see Proposition 1.5 below)
$$K(G/B)\cong R(T)/\cI,$$
where $\cI$ is the ideal generated by 
$\{ f\in R(T)^W \ |\ f-\varepsilon(f)=0\}$.  Equivalently,
$K(G/B)\cong R(T)\otimes_{R(G)}\ZZ$,
where $R(G)$ acts on $\ZZ$ by $[V]\cdot 1 = \dim(V)$, if $V$ is 
a $G$-module.

\vfill\eject 

\medskip
\subsection{$K(G/P)$ for a parabolic subgroup $P$}

A similar setup works when $B$ is replaced by any parabolic subgroup $P$
containing $B$. The coset space $G/P$ is a projective variety and the Bruhat
decomposition takes the form $G = \bigcup_{\bar w\in W/W_P} B\bar wP$
where $W_P$ is the subgroup of $W$ given by
$W_P=\langle s_i\ |\ \gg_{-\alpha_i}\in \gp\rangle$,
where $\gp$ is the Lie algebra of $P$.
The Schubert varieties $X_{\bar w}$ are the closures of the 
Schubert cells $X_{\bar w}^\circ=B\bar wP$ in $G/P$.
$$\hbox{$K(G/P)$ is a free $\ZZ$-module with basis 
$\{ [\cO_{X_{\bar w}}]\ |\ \bar w\in W/W_P\}$. }$$

Write $P=LU$ where $U$ is the unipotent radical of $P$ and $L$ is a Levi
subgroup.  The Weyl group of $L$ is $W_P$ and
$$R(L) \cong R(T)^{W_P}$$
is the subring of $W_P$-symmetric functions in $R(T)$.
The same construction as in 
(1.2) with $B$ replaced by $P$ and $T$ replaced by $L$ gives
a  ring homomorphism 
$$\matrix{
\phi_P\colon &R(L) &\longrightarrow &K(G/P) \cr
&V &\longmapsto &(G\times_P V \mapright{\pi} G/P). \cr}
\formula$$ 

\prop  Let $G$ be a connected simply connected semisimple Lie group and
let $T$ be a maximal torus of $G$.  
Let $P$ be a parabolic subgroup of $G$ with Levi decomposition $P=LU$
and let $W_P$ be the Weyl group of $L$.  Then
$$K(G/P)\cong {R(T)^{W_P}\over \cI_P},$$
where $\cI_P$ is the ideal generated by 
$\{ f\in R(T)^{W_P} \ |\ f-\varepsilon(f)=0\}$ and $\varepsilon\colon R(T)\to
\ZZ$ is the map given by $\varepsilon(e^\lambda)=1$ for $\lambda\in \Lambda$. 
Equivalently,
$K(G/P)=R(L)\otimes_{R(G)}\ZZ$.
\pf 
Let $K_{vb}(G/P)$ be the Grothendieck group of $C^\infty$ vector bundles on
$G/P$ and let $\eta\colon K(G/P)\to K_{vb}(G/P)$ be the map which assigns to a
locally free sheaf its underlying vector bundle.
Let $\tilde \phi_P$ be the  composition 
$$\tilde\phi_P\colon R(L)\mapright{\phi_P} K(G/P)\mapright{\eta} K_{vb}(G/P).$$
Since $\pi_1(G)=0$, $\pi_1(P)\cong \pi_2(G/P)$, which is free abelian by the
Bruhat decomposition. Since the unipotent radical $U$ of $P$ is contractible the
projection 
$f\colon P\to P/U\cong L$ is a homotopy equivalence.  
Thus $\pi_1(L)$ is free abelian and  
we may apply the results of [P] to conclude that $\tilde \phi_P\colon R(L)\to
K_{vb}(G/P)$ is surjective, $R(L)$ is projective over $R(G)$ with rank
$|W/W_P|$ and $K_{vb}(G/P)$ is a free $\ZZ$-module of the same rank.
(Note:  The results of [P] can be applied since $G$ and $L$ are the
complexifications of compact groups.)

Since $\tilde \phi_P$ is surjective the map
$\eta$ is also surjective.   Then, since $K(G/P)$ and $K_{vb}(G/P)$ are both free
$\ZZ$-modules of rank
$|W/W_P|$, it follows that $\eta$ must be an isomorphism.  
This means two things: (1) that we can identify $K(G/P)$ and $K_{vb}(G/P)$,
and (2) that $\phi_P$ is surjective.

The kernel $\cI_P$ of $\phi_P$ is identified by using (1.3).
\endpf

\vfill\eject

\subsection{Transfer from $K(G/B)$ to $K(G/P)$}

Although we will work primarily with $K(G/B)$ it is standard to transfer 
results from $K(G/B)$ to results on $K(G/P)$.  This can be accomplished
with the following proposition.  The proof will be given in section 5.

\prop If $f\colon G/B\to G/P$ is the natural projection then the induced map
$f^!\colon K(G/P)\to K(G/B)$ is an injection.
This map is given explicitly by
$$f^!([\cO_{X_{\bar w}}])=[\cO_{X_v}],$$
where $v\in W$ is the unique element of longest length in the coset
$\bar w=vW_P$.
\endprop

\section 2. The class $[\cO_{P/B}]$ in $K(G/B)$ 

In this section we give an expression for the class
$[\cO_{P/B}]\in K(G/B)$ as an element of $R(T)/\cI$.  This is done by first
finding a formula for
$[\cO_{P/P}]$ in $K(G/P)$ and then using the projection
$f\colon G/B\to G/P$ to pull back this formula to
$K(G/B)$. The formula for $[\cO_{P/P}]$ in $K(G/P)$ is obtained by using a Koszul
resolution on a vector field with simple zeros at the points 
$\{\bar w_iP\ |\ \bar w_i\in W/W_P\}$.  
This reduces the computation to determining
$\Lambda_{-1}(T^*(G/P))$ and this can be done since we understand the structure
of $T^*(G/P)=(\gg/\gp)^*$ as a $B$-module (under the adjoint action).
Although these formulae for $[\cO_{P/B}]$ are useful for specific computations
they are not needed for the proof of our main result, Theorem 4.3.
 
\thm  Let $P\supseteq B$ be a parabolic subgroup of $G$ and let 
$w$ be the longest element of the corresponding parabolic subgroup $W_P$ of $W$. 
In $K(G/B)$
$$[\cO_{X_w}]=\hbox{$|W_P|\over|W|$}
\prod_{\gg_{-\alpha}\not\in \gp} (1-e^{-\alpha}),$$
where the product is over all positive roots $\alpha$ such that 
$\gg_{-\alpha}\not\in \gp$.
\pf
Let $\gh$ denote the complex Lie algebra of
the maximal torus $T\subseteq G$ and let $H\in \gh$ be a {\it regular} element,
i.e. the $W$ action on $H$ has trivial stabilizer.  The
one-parameter group $exp(zH)$, $z\in \CC$, of $G$ induces a flow on $G/P$ (by
left translation) whose fixed points are the points in the set 
$$Z=\{w_iP\in G/P\},$$ 
where $w_i$ run over a set of coset representatives of $W/W_P$.  It follows that
the zeros of the associated vector field $v(H)$ are the same points and a local
calculation shows that they are simple.  This construction of vector fields
$v\colon G/P\to T(G/P)$ whose zeros are isolated and simple is essentially due to
A. Weil [W].

Since the zero set $Z$ of the vector field $v(H)$ is a smooth subvariety of
codimension equal to the fibre dimension of $T(G/P)$, the vector field $v(H)$
gives rise to a Koszul resolution of $\cO_Z$ (see [CG] \S 5.4)
$$\hbox{$\cdots
\mapright{i_v} \cO_{G/P}(\bigwedge^{2}(T^*(G/P)))
\mapright{i_v} \cO_{G/P}(\bigwedge^{1}(T^*(G/P)))
\longrightarrow\cO_{G/P}\longrightarrow\cO_{Z}
\longrightarrow 0,$}$$
where 
$i_v$ denotes interior product with $v(H)$.
Hence, in $K(G/P)$ we have 
$$[\cO_Z]=\Lambda_{-1}(T^*(G/P))
\qquad\hbox{where}\qquad
\Lambda_{-1}(T^*(G/P)) = \sum_i \hbox{$(-1)^i \bigwedge^i(T^*(G/P))$}.$$

For any two points $p,q\in G/P$ 
there is a $g\in G$ so that
$gp=q$.  Since $G$ acts trivially on $K(G/P)$,
$$[{\cal O}_q] = g^![{\cal O}_q] = [g^*{\cal O}_q] = [{\cal O}_p].\formula$$
Since the points $\bar w_iP$ of $Z$ are simple, 
$$\cO_Z=
\bigoplus_{i=1}^{|W/W_P|} \cO_{\bar w_iP}
\qquad\hbox{and so}\qquad 
[\cO_Z]=\sum_{i=1}^{|W/W_P|}[\cO_{\bar w_iP}]=|W/W_P|[\cO_P],$$ 
by (2.2).  Since $K(G/P)$ is a free 
$\ZZ$-module it follows that $[\cO_Z]$ is divisible
by $|W/W_P|$ and we get
$$[\cO_P]={|W_P|\over |W|} \Lambda_{-1}(T^*(G/P)).\formula$$

Let us compute the pull back
$f^!(|W/W_P|[\cO_P])=f^!(\Lambda_{-1}(T^*(G/P)))\in K(G/B)$ for the projection
$f\colon G/B\to G/P$.  The bundle
$T^*(G/P)$ is the homogeneous vector bundle over $G/P$ associated to the
$P$-module $(\gg/\gp)^*$, where $P$ acts on $\gg/\gp$ by the adjoint action.  Then
$f^!(T^*(G/P))$ is the vector bundle over
$G/B$ associated to the $B$-module $(\gg/\gp)^*$, 
where we regard $(\gg/\gp)^*$ as a
$B$-module by restriction.  By Lie's theorem,
$(\gg/\gp)^*$ admits an $B$-module filtration such that the unipotent
radical of $B$ acts trivially on the associated graded module $\gr_F(\gg/\gp)^*$. 
Hence
$$
\gr_F(\gg/\gp)^*=\sum_{ \gg_{-\alpha}\not\in \gp} \gg_{\alpha} $$
is a sum of weight spaces as an $\ad(B)$-module.
Since a filtered object and its associated graded define the same element in a
Grothendieck ring we have 
$$f^![T^*(G/P)]=\sum_{\gg_{-\alpha}\not\in \gp} e^{-\alpha}$$
in $K(G/B)$.  From this equation we get the formula for
$$f^!(\Lambda_{-1}(T^*(G/P))=\prod_{\gg_{-\alpha}\not\in \gp} (1-e^{-\alpha}).
\formula$$  
The theorem follows from (2.3), (2.4) and Proposition 1.6 since
$$f^!([\cO_P])=[\cO_{f^{-1}(P)}]=[\cO_{P/B}]$$
and $P/B=X_w$ for the longest element $w$ of $W_P\subseteq
W$.
\endpf

\cor In $K(G/B)$
\smallskip
\itemitem{} $\displaystyle{[\cO_{X_1}] = {1\over |W|}
\prod_{\alpha>0} (1-e^{-\alpha}) }$,
\smallskip
\itemitem{} $\displaystyle{
[\cO_{X_{s_i}}] = {2\over |W|} \prod_{\alpha>0\atop \alpha\ne \alpha_i}
(1-e^{-\alpha}), }$
\quad for each simple reflection $s_i$, $1\le i\le n$, and
\smallskip
\itemitem{} $\displaystyle{[\cO_{X_{w_0}}] = 1, }$ \quad
for the longest element $w_0$ in $W$.
\pf
The first and last formulas are Theorem 2.1 in the cases 
$P=B$ and $P=G$ respectively and the middle formula is the case when $P=P_i$ is
the minimal parabolic subgroup  whose Lie algebra $\gp_i$ is generated by $\gb$
and the negative root space $\gg_{-\alpha_i}$.
\endpf

\noindent
{\bf Remarks.}  
\smallskip\noindent
{\bf 1.}  In optimal cases such as $G=SL(n,\CC)$
one can use various tautological bundles on $SL(n,\CC)/P$
to construct resolutions of $\cO_p$ directly, and hence obtain
formulae for $[\cO_p]$ which are ``denominator free''.  One example is obtained
from the tautological
$k$-plane bundle over the Grassmannian of $k$-planes in $\CC^n$: 
$E_k\longrightarrow \GG(k,\CC^n)$.  Every (homogeneous) linear function on
$\CC^n$ defines an algebraic section of 
$E_k^*$.  Hence by choosing $(n-k)$ linearly independent such functions,
we can define a section $\sigma\colon \GG(k,\CC^n)\to 
\bigoplus_{n-k} E_k^*$ whose unique zero is the point $p$ corresponding
to the common kernel of the linear functions.  Since $\sigma$ is clearly
regular, $[\cO_p] = \sum (-1)^i[\bigwedge^i(\bigoplus_{n-k}E_k)]$
in $K(\GG(k,\CC^n))$.  Other examples can be found in [FL].

\medskip\noindent
{\bf 2.}  In contrast with the previous remark, it seems difficult to find
``denominator  free'' formulae for $[\cO_p]$ in general.  A comparison with
cohomology will be helpful.  For $\FF(\CC^n)=SL(n,\CC)/B$ a generator of the top
cohomology is given by $x_1^{n-1}x_2^{n-2}\cdots x_{n-1}$, where 
$x_j\in H^2(\FF(\CC^n);\ZZ)$ form a suitable basis.  For general
$G/B$ the only uniform expressions for a generator in the top degree
all involve denominators.  For example, one such is 
${1\over|W|}\prod_{\alpha>0}\alpha$.  Indeed, if $H^*(G;\ZZ)$ has torsion
then no integral polynomial in a basis for $H^2(G/B;\ZZ)$ will give
a generator in the top degree.

\section 3.  Push-pull operators in K-theory
 
For a positive root $\alpha$, let $s_\alpha\in W$ be the corresponding 
reflection and define operators $L_\alpha\colon R(T)\to R(T)$ and
by $T_\alpha\colon R(T)\to R(T)$ by
$$L_\alpha(x)={x-s_\alpha x\over 1-e^{-\alpha}}
\qquad\hbox{and}\qquad
T_\alpha(x) = e^{-\rho}L_\alpha(e^\rho x)
={e^\alpha x-s_\alpha x\over e^\alpha - 1},$$ 
respectively, where $\rho={1\over2}\sum_{\alpha>0} \alpha$.
In this section we will show that there is an inductive formula
for the classes $[\cO_{X_w}]$ in terms of the
operators $T_\alpha$ and the class $[\cO_{X_1}]$,
which was determined in Corollary 2.5.

The operators $T_\alpha$ and $L_\alpha$ have been in the literature for some 
time, see for example, [D, \S 5], [KK], [FL]. 
Let $x,y\in R(T)$.  Short direct calculations using the definitions
establish the following identities:
\global\advance\resultno by 1
\medskip
\itemitem{(\the\sectno.\the\resultno a)}  
$T_\alpha(xy)=xT_\alpha(y)$, \quad if $s_\alpha x=x$,
\medskip
\itemitem{(\the\sectno.\the\resultno b)}  
$s_\alpha T_\alpha x = T_\alpha x$,
\medskip
\itemitem{(\the\sectno.\the\resultno c)}  
$T_\alpha T_\alpha x = T_\alpha x$,
\medskip
\itemitem{(\the\sectno.\the\resultno d)}  
$\displaystyle{ 
e^\lambda T_\alpha x = 
\left( T_\alpha e^{s_\alpha\lambda} + 
{e^\lambda-e^{s_\alpha\lambda}\over 1-e^{-\alpha}}\right)x. }$
\medskip\noindent
Because of (a), $T_\alpha$ is a map of $R(G)$-modules and so it
descends to an operator on $K(G/B)$ which we shall denote by the same symbol.
Moreover, the induced operator on $K(G/B)$ satisfies (3.1a-d). 

Let $\alpha$ be a simple root and let $P_\alpha$ be the minimal parabolic 
subgroup whose Lie algebra
${\goth p}_\alpha$ is generated by ${\goth b}$ and the negative root space 
${\goth g}_{-\alpha}$.  Since $P_\alpha/B\cong \PP_1$, the natural projection 
$$f_\alpha\colon G/B\to G/P_\alpha\formula$$
is a $\PP_1$-bundle.

\prop  Let $\alpha$ be a simple root.
For every $x\in K(G/B)$, 
$$(f_\alpha)^!\circ(f_\alpha)_!(x)=T_\alpha(x).$$
\endprop

This result is proved in [KK, Prop. 4.11].  In section 5
we shall see that the Grothendieck-Riemann-Roch theorem implies that this fact
is equivalent to the corresponding fact in cohomology.
This alternate point of view has the advantage that it illustrates why the
operators $T_\alpha$ are the K-theoretic analogues of the BGG operators
$\partial_\alpha$ (see [BGG] and [D]).  The proof of following proposition is a
generalization of the argument in [FL, p. 728].  Kostant and Kumar [KK, Lemma
4.12] have also proved the same result.

\prop Let $s_\alpha\in W$ be the simple reflection corresponding to a
simple root $\alpha$. Given a Schubert variety $X_w\subseteq G/B$,
$$(f_\alpha)^!\circ(f_\alpha)_!([\cO_{X_w}])
=\cases{
[\cO_{X_{ws_\alpha}}], &if $\ell(ws_\alpha)>\ell(w)$, \cr
[\cO_{X_w}], &if $\ell(ws_\alpha)<\ell(w)$. \cr}
$$
\pf  The main idea of the proof is
$$(f_\alpha)^!\circ (f_\alpha)_!([\cO_{X_w}])
=(f_\alpha)^!([\cO_{f_\alpha(X_w)}])
=[\cO_{f_\alpha^{-1}(f_\alpha(X_w))}].$$
One only has to justify the equalities and identify $f_\alpha(X_w)$ and
$f_\alpha^{-1}(f_\alpha(X_w))$.

For $w\in W$ let $\bar w=\{w,ws_\alpha\}$. 
It is convenient to relabel the elements of the set $\{w,ws_\alpha\}$ as
$w'$ and $w''$ where by fiat $\ell(w'')=\ell(w')+1$.  
Analyzing the Bruhat decomposition of $X_w$ in (1.1) we get
$$f_\alpha(X_{w'})=f_\alpha(X_{w''})=X_{\bar w}\qquad\hbox{and}\qquad 
f_\alpha^{-1}(X_{\bar w}) = X_{w''}.\formula$$
Since $f_\alpha\colon X_{w'}^\circ\to X_{\bar w}^\circ$ is an isomorphism of
varieties $f_\alpha\colon X_{w'}\to X_{\bar w}$ is birational.  This
combined with the (deep) fact that Schubert varieties have rational
singularities (see the survey [Ra] and the references there) implies that 
\smallskip\noindent
\itemitem{(a)}  $(f_\alpha)_*(\cO_{X_{w'}})=\cO_{X_{\bar w}}$,
and 
$R^q(f_\alpha)_*(\cO_{X_{w'}})=0$, for $q>0$.
\smallskip\noindent
From the Bruhat decomposition one sees that  
$f_\alpha\colon X_{w''}\to X_{\bar w}$ is the restriction of the ambient
$\PP_1$-bundle $f_\alpha\colon G/B\to G/P_\alpha$.
Thus 
\smallskip\noindent
\itemitem{(b)} $(f_\alpha)_*(\cO_{X_{w''}})=\cO_{X_{\bar w}}$ and
$R^q(f_\alpha)_*(\cO_{X_{w''}})=0$, for $q>0$.
\smallskip\noindent
Finally, from (3.5) we have
\smallskip\noindent
\itemitem{(c)} $(f_\alpha)^*(\cO_{X_{\bar w}})=\cO_{X_{w''}}$.
\smallskip\noindent
\smallskip\noindent
Statements (a) and (b) imply that 
$(f_\alpha)_!([\cO_{X_w}])=[\cO_{f_\alpha(X_w)}]$
and (c) implies that $(f_\alpha)^!([\cO_{X_{\bar w}}])=[\cO_{X_{w''}}]$.
\endpf

\cor  For each simple root $\alpha_i$ let $T_i=T_{\alpha_i}$.
Let $w=s_{i_1}\cdots s_{i_p}$ be a reduced expression for
$w$ and define $T_w = T_{i_1}\cdots T_{i_p}$.  Then $T_w$ is
independent of the choice of the reduced expression of $w$ and
$$T_{w^{-1}}[\cO_{X_1}] = [\cO_{X_w}].$$
\pf
The formula in the statement follows from Propositions 3.3 and 3.4.
These two Propositions, combined with formula (3.1c)
also show that the action of $T_w$ on the elements of the basis 
$\{[\cO_v]\ |\ v\in W\}$ of $K(G/B)$ is independent
of the choice of the reduced word for $w$.  By
Proposition 1.5, $K(G/B)$ is a free $R(G)$-module and thus
it follows from (3.1c) that, as an operator on
$R(T)$, $T_w$ is independent of the reduced word for $w$. 
\endpf

\vfill\eject

\section 4.  The Pieri-Chevalley formula

In this section we shall inductively apply formula (3.1d) to obtain
an expansion of the product $e^\lambda[\cO_{X_w}]$ in $K(G/B)$ in terms of the
basis $\{[\cO_{X_v}]\ |\ v\in W\}$.  We use the path model of P. Littelmann
to keep track of the combinatorics involved in iterating formula (3.1d).
 
\subsection{The path model}

Let $\Lambda=\sum_i \ZZ\omega_i$ be the weight lattice and let $\gh^*=\sum_i
\RR\omega_i$.  A {\it path} in $\gh^*$ is a piecewise linear map 
$\pi\colon [0,1]\to \gh^*$ such that $\pi(0)=0$.  Let $\pi$ be a path,
let $\alpha$ be a simple root and let $h_\alpha\colon [0,1]\to \RR$ be the
function given by
$$\matrix{
h_\alpha\colon &[0,1] &\longrightarrow &\RR \cr
&t &\longmapsto &\langle\pi(t),\alpha^\vee\rangle. \cr}$$ 
At $t$
this function gives the position of $\pi(t)$ in the $\alpha$-direction. Let
$m_\alpha$ be the minimal value of $h_\alpha$ and define functions
$l\colon [0,1]\to [0,1]$ and $r\colon [0,1]\to [0,1]$ by
$$l(t)={\rm min}\{1,h_\alpha(s)-m_\alpha\ | t\le s\le 1\},
\qquad
r(t)=1-{\rm min}\{1,h_\alpha(s)-m_\alpha\ | 0\le s\le t\}.$$
The {\it root operators} (see [L3] Definitions 2.1 and 2.2) are 
operators on the paths given by
$$\eqalign{
e_\alpha\pi &= \cases{
t\mapsto \pi(t)+r(t)\alpha, &if $r(0)=0$, \cr
0, &otherwise, \cr} \qquad\hbox{and} \cr
\cr
f_\alpha\pi &= \cases{
t\mapsto \pi(t)-l(t)\alpha, &if $l(1)=1$, \cr
0, &otherwise, \cr}
}
$$
where we use $0$ to denote the ``null path''.  

Fix a dominant weight $\lambda\in \Lambda$.  
Let $\pi_\lambda$ be the path given by
$\pi_\lambda(t)=t\lambda$, $0\le t\le 1$, and let
$$\cT^\lambda = \{ f_{i_1}f_{i_2}\cdots f_{i_l}\pi_\lambda \}$$ 
be the set of all paths obtained by applying sequences of
root operators $f_i=f_{\alpha_i}$, $1\le i\le n$ to $\pi_\lambda$.  This is the
set of  Lakshmibai-Seshadri paths of shape $\lambda$. P. Littelmann [L1] has
shown that this set of paths is finite and can be characterized in terms of an
integrality condition.  We shall not need this alternative characterization. 

Let $W_\lambda$ be the stabilizer of $\lambda$.  The cosets in $W/W_\lambda$ are
partially ordered by the Bruhat-Chevalley order.  Use a pair of sequences
$$\matrix{
\hfill\vec \tau&=(\tau_1>\tau_2>\cdots >\tau_\ell),\hfill
&\tau_i\in W/W_\lambda, \qquad \hbox{and}\hfill\cr
\hfill\vec a&=(0=a_0<a_1<a_2<\cdots<a_\ell=1),\hfill
&\qquad a_i\in \QQ,\hfill\cr
}$$
to encode the path $\pi\colon [0,1]\to \gh^*$ given by
$$\pi(t)=(t-a_{j-1})\tau_j\lambda+\sum_{i=1}^{j-1} (a_i-a_{i-1})\tau_i\lambda,
\qquad\hbox{for $a_{j-1}\le t\le a_j$.}
$$
We shall write $\pi=(\vec \tau,\vec a)$.  Every path $\pi\in \cT^\lambda$ is
of this form.
Littelmann introduced this set of paths $\cT^\lambda$
as a model for the Weyl character formula.  He proved that
$$\sum_{\eta\in \cT^\lambda} e^{\eta(1)}
={\sum_{w\in W} \varepsilon(w)e^{w(\lambda+\rho)}
\over \sum_{w\in W} \varepsilon(w)e^{w\rho} },
$$
where $\rho={1\over2}\sum_{\alpha>0}\alpha$ is the half-sum of the positive roots.

\vfill\eject

\subsection{Application of the path model}

Fix a dominant weight $\lambda\in \Lambda$ and let $\pi=(\vec\tau,\vec a)
=((\tau_1>\cdots>\tau_r),(a_0<\cdots<a_r))\in \cT^\lambda$.  The {\it
initial direction} of $\pi$ is $\iota(\pi)=\tau_1$. 
Fix $w\in W$, let $\bar w=wW_\lambda\in W/W_\lambda$ and define
$$\cT^\lambda_w = \{ \pi\in \cT^\lambda \ |\ \iota(\pi)\le \bar w\}.$$
Let $\pi=(\vec\tau,\vec a)\in \cT^\lambda_w$. 
A {\it maximal lift of $\vec\tau$ with respect to $w$} is a choice of 
representatives $t_i\in W$ of the cosets $\tau_i$ such that 
$w\ge t_1>\cdots>t_r$ and each $t_i$ is maximal in Bruhat order such that
$t_{i-1}>t_i$.   The {\it final direction} of $\pi$ with respect to $w$ is
$$v(\pi,w)=t_r,$$
where
$w\ge t_1>\cdots>t_r$ is a maximal lift of $\tau_1>\ldots>\tau_r$
with respect to $w$.

For each  $\pi\in \cT^\lambda$ such that $e_\alpha(\pi)=0$
the {\it $\alpha$-string of $\pi$} is the set of paths
$$S_\alpha(\pi) = \{ \pi, f_\alpha\pi, \ldots, f^m_\alpha\pi\},$$
where $m$ is maximal such that $f^m_\alpha\pi\ne 0$.  We have:
\smallskip
\itemitem{(a)} 
If $f_\alpha^j\pi\ne 0$ then $(f_\alpha^j \pi)(1)=\pi(1)-j\alpha$. 
\smallskip
\itemitem{(b)} $\iota(f_\alpha^j\pi)=s_\alpha\iota(\pi)$ for all $1\le j\le m$.
\smallskip
\itemitem{(c)} If $S_\alpha(\pi)\subseteq \cT^\lambda_w$ then
$v(f_\alpha^m\pi,w)=s_\alpha v(\pi,w)$  and $v(f_\alpha^j \pi,w)=v(\pi,w)$
for
$1\le j<m$.
\smallskip\noindent
Statement (a) is [L2] Lemma 2.1a, statement (b) is [L1] Lemma 5.3b, and 
statement (c) follows from [L1] Lemma 5.3c and [L2] Lemma 2.1e.  All of these
facts are really coming from the explicit form of the action of the root
operators on the Lakshmibai-Seshadri paths which is given in [L1] Proposition
4.2.  The consequence of (a) and (c) is that 
$$\sum_{\eta\in S_\alpha(\pi)} T_{v(\eta,w)^{-1}} e^{\eta(1)}
=T_{v(\pi,w)^{-1}}\left(T_\alpha e^{s_\alpha\pi(1)} +
{e^{\pi(1)}-e^{s_\alpha\pi(1)}\over 1-e^{-\alpha}}\right)
=T_{v(\pi,w)^{-1}} e^{\pi(1)} T_\alpha.
$$

Let $w=s_\alpha w'$ where $\ell(w)=\ell(w')+1$. 
Let $\pi\in \cT^\lambda_w$ be such that $e_\alpha(\pi)=0$.
It follows from (a) that 
$$\hbox{$S_\alpha(\pi)\subseteq \cT^\lambda_{w}$,\quad and \qquad
either $S_\alpha(\pi)\cap \cT^\lambda_{w'}=\{\pi\}$ or
$S_\alpha(\pi)\subseteq \cT^\lambda_{w'}$.}$$
Suppose that $w\ge t_1>\cdots>t_r$ and $w'\ge t_1'>\cdots>t_r'$ are maximal
lifts of $\pi$ with respect to $w$ and $w'$ respectively.
\smallskip
\item{} If $m>0$ then $t_1$ is not divisible by $s_\alpha$.  It follows that
$t_1'=t_1$ and thus that $t_r=t_r'$.
\smallskip
\item{}
If $m=0$ then all the $t_i$ are divisible by $s_\alpha$ and it follows that
$t_r=s_\alpha t_r'$. 
\smallskip\noindent
Thus $v(\pi,w)=v(\pi,w')$ if $m>0$ and $v(\pi,w)=v(\pi,w')s_\alpha$ if
$m=0$. We conclude that 
$$\eqalign{
T_{v(\pi,w')^{-1}} e^{\pi(1)} T_\alpha 
&=T_{v(\pi,w')^{-1}}\left(T_\alpha e^{s_\alpha\pi(1)} +
{e^{\pi(1)}-e^{s_\alpha\pi(1)}\over 1-e^{-\alpha}}\right)  \cr
&=\sum_{\eta\in S_\alpha(\pi)} T_{v(\eta,w)^{-1}} e^{\eta(1)}. \cr
}
\formula$$

\thm  Let $\lambda$ be a dominant integral weight and let $w\in W$.  Then
$$e^\lambda T_{w^{-1}} = 
\sum_{\eta\in \cT^\lambda_w} T_{v(\eta,w)^{-1}} e^{\eta(1)}
$$
as operators on $R(T)$.
\pf
The proof is by induction on $\ell(w)$.  The base case $\ell(w)=1$ is
formula (3.1d).   Let
$w=s_\alpha w'$ with $\ell(w)=\ell(w')+1$.  Then 
$$\eqalign{
e^\lambda T_{w^{-1}}
&=e^\lambda T_{(w')^{-1}}T_\alpha \cr
&=\left(\sum_{\eta\in \cT^\lambda_{w'}} 
T_{v(\eta,w')^{-1}}e^{\eta(1)}\right) T_\alpha
\qquad\qquad\qquad\hbox{(by induction)} 
\cr 
&= \sum_{\pi\in \cT^\lambda_w\atop e_\alpha(\pi)=0} 
\left(\sum_{S_\alpha(\pi)\subseteq \cT^\lambda_{w'}}
T_{v(\pi,w')^{-1}}e^{\pi(1)}T_\alpha 
+
\sum_{S_\alpha(\pi)\cap \cT^\lambda_{w'}=\{\pi\} }
T_{v(\pi,w')^{-1}}e^{\eta(1)} \right)T_\alpha \cr
&= \sum_{\pi\in \cT^\lambda_w\atop e_\alpha(\pi)=0} 
\left(\sum_{S_\alpha(\pi)\subseteq \cT^\lambda_{w'}}
T_{v(\pi,w')^{-1}}e^{\pi(1)}T_\alpha +
\sum_{S_\alpha(\pi)\cap \cT^\lambda_{w'}=\{\pi\} }
T_{v(\pi,w')^{-1}}e^{\eta(1)}T_\alpha \right) \cr
&= \sum_{\eta\in \cT^\lambda_w} T_{v(\eta,w)^{-1}} e^{\eta(1)}
\qquad\qquad\qquad\qquad\hbox{(by (4.1))}.\qquad\hbox{\qed} \cr
}$$

\thm  Let $\lambda$ be a dominant integral weight. In $K(G/B)$
$$e^\lambda [\cO_{X_w}] = 
\sum_{\eta\in \cT^\lambda_w} [\cO_{X_{v(\eta,w)}}] 
$$
\vskip-.2in
\pf
Since the sheaf $\cO_{X_1}$ is supported on the single
point $X_1\in G/B$, any product $[\cF][\cO_{X_1}]$, where $\cF$ is a vector
bundle on $G/B$ is the class of a bundle supported on the single point $X_1$.  
More precisely, $[\cF][\cO_{X_1}]=\hbox{rk}(\cF)[\cO_{X_1}]$.  Thus, since
$e^\lambda$ is the class of a line bundle we have
$e^\lambda[\cO_{X_1}]=[\cO_{X_1}]$.   By Corollary 3.6, 
$[\cO_{X_w}]=T_{w^{-1}}[\cO_{X_1}]$ and so the result follows from Theorem
4.2.
\endpf

The following example illustrates how one computes the product
$e^{\omega_2}[\cO_{X_{s_1s_2s_1s_2}}]$ in $K(G/B)$ for the group $G$
of type $G_2$.  In this case $\lambda=\omega_2$, $w^{-1}=s_2s_1s_2s_1$
and the starting path
$\pi_\lambda$ is the straight line path from the origin to the point
$\omega_2$.  The paths in the set $\cT^{\omega_2}_{s_2s_1s_2s_1}$ are
the paths in the following diagrams.
$$
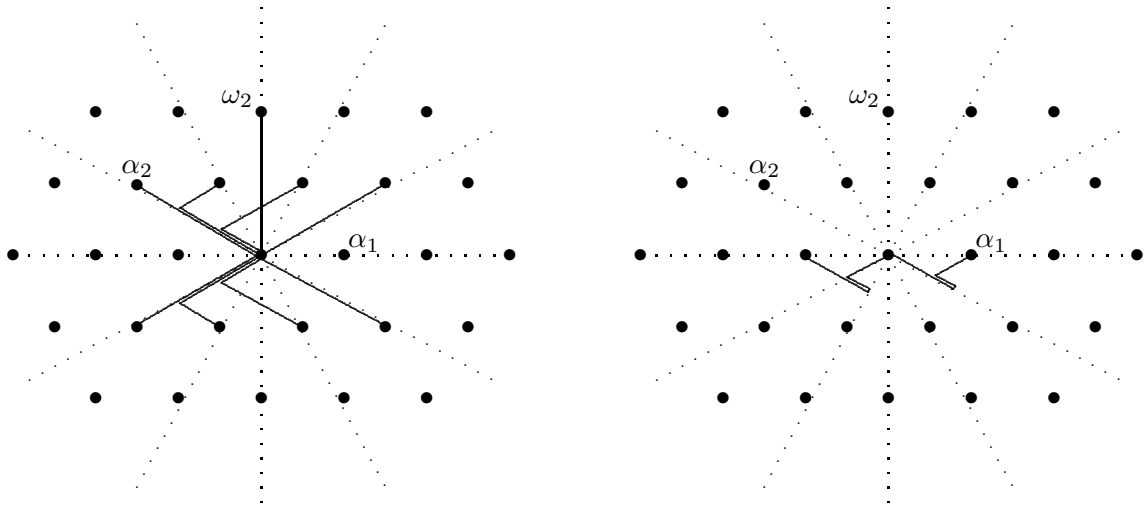
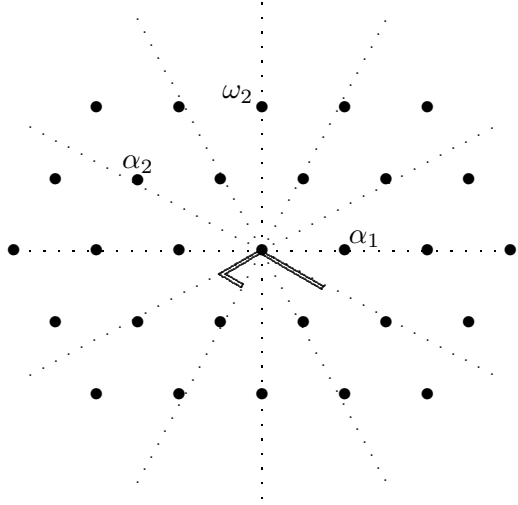
\beginpicture
\setcoordinatesystem units <1cm,1cm>         
\setplotarea x from -3 to 3, y from -3 to 3    
\put{$\omega_2$}[r] at -0.1 1.9
\put{$\alpha_2$} at -1.5 1.05
\put{$\alpha_1$}[l] at 1.05 0.15
\put{$\bullet$} at 0 0
\put{$\bullet$} at 1 0
\put{$\bullet$} at 2 0
\put{$\bullet$} at 3 0
\put{$\bullet$} at -1 0
\put{$\bullet$} at -2 0
\put{$\bullet$} at -3 0
\put{$\bullet$} at 0.5  0.8660
\put{$\bullet$} at 1.5  0.8660
\put{$\bullet$} at 2.5  0.8660
\put{$\bullet$} at -0.5  0.8660
\put{$\bullet$} at -1.5  0.8460
\put{$\bullet$} at -2.5  0.8660
\put{$\bullet$} at 0  1.732
\put{$\bullet$} at 1  1.732
\put{$\bullet$} at 2  1.732
\put{$\bullet$} at -1  1.732
\put{$\bullet$} at -2  1.732
\put{$\bullet$} at 0.5  -0.8660
\put{$\bullet$} at 1.5  -0.8660
\put{$\bullet$} at 2.5  -0.8660
\put{$\bullet$} at -0.5  -0.8660
\put{$\bullet$} at -1.5  -0.8660
\put{$\bullet$} at -2.5  -0.8660
\put{$\bullet$} at 0  -1.732
\put{$\bullet$} at 1  -1.732
\put{$\bullet$} at 2  -1.732
\put{$\bullet$} at -1  -1.732
\put{$\bullet$} at -2  -1.732
\linethickness=0.5pt                          
\putrule from 0 1.732 to 0 0          
\plot 0 0 1.5  0.8660 /    
\plot 0 0.05 -0.48  0.3187 /    
\plot  -0.48  0.3187 0.5 0.8660 /    %
\plot 0 -0.00 -.98  0.57735 /    
\plot  -.98  0.57735 -0.5 0.8660 /    %
\plot 0 -0.04 -1.5  0.8460 /    
\plot 0 0.03 -1.5  -0.8360 /    
\plot 0 0.0 -.98  -0.57735  /    
\plot  -.98  -0.57735  -0.5 -0.8660 /    %
\plot 0 -0.04 -0.48  -0.3287 /    
\plot  -0.48  -0.3287  0.5 -0.8660 /    %
\plot 0 -0.04 1.5  -0.8460 /    
\setdots
\putrule from 0 3 to 0 -3        
\plot -1.5 -2.798   1.5 2.798 /  
\putrule from 3 0  to -3 0       
\plot -2.798 -1.5  2.798 1.5 /   
\plot -2.798 1.5  2.798 -1.5 /   
\plot -1.5 2.798  1.5 -2.798 /   
\endpicture
\qquad\qquad
\beginpicture
\setcoordinatesystem units <1cm,1cm>         
\setplotarea x from -3 to 3, y from -3 to 3    
\put{$\omega_2$}[r] at -0.1 1.9
\put{$\alpha_2$} at -1.5 1.05
\put{$\alpha_1$}[l] at 1.05 0.15
\put{$\bullet$} at 0 0
\put{$\bullet$} at 1 0
\put{$\bullet$} at 2 0
\put{$\bullet$} at 3 0
\put{$\bullet$} at -1 0
\put{$\bullet$} at -2 0
\put{$\bullet$} at -3 0
\put{$\bullet$} at 0.5  0.8660
\put{$\bullet$} at 1.5  0.8660
\put{$\bullet$} at 2.5  0.8660
\put{$\bullet$} at -0.5  0.8660
\put{$\bullet$} at -1.5  0.8460
\put{$\bullet$} at -2.5  0.8660
\put{$\bullet$} at 0  1.732
\put{$\bullet$} at 1  1.732
\put{$\bullet$} at 2  1.732
\put{$\bullet$} at -1  1.732
\put{$\bullet$} at -2  1.732
\put{$\bullet$} at 0.5  -0.8660
\put{$\bullet$} at 1.5  -0.8660
\put{$\bullet$} at 2.5  -0.8660
\put{$\bullet$} at -0.5  -0.8660
\put{$\bullet$} at -1.5  -0.8660
\put{$\bullet$} at -2.5  -0.8660
\put{$\bullet$} at 0  -1.732
\put{$\bullet$} at 1  -1.732
\put{$\bullet$} at 2  -1.732
\put{$\bullet$} at -1  -1.732
\put{$\bullet$} at -2  -1.732
\linethickness=0.5pt                          
\plot 0 0.04 0.78  -0.40 /    
\plot 0.78  -0.40  0.81  -0.37 /    %
\plot 0.81  -0.37  0.57  -0.24 /    %
\plot  0.57  -0.24  1 0 /    %
\plot 0 0 -0.50 -0.26  /  
\plot -0.50 -0.26 -0.23 -0.40 /  %
\plot -0.23 -0.40 -0.25 -0.44 /  %
\plot -0.25 -0.44 -1.04 0 /  %
\setdots
\putrule from 0 3 to 0 -3        
\plot -1.5 -2.798   1.5 2.798 /  
\putrule from 3 0  to -3 0       
\plot -2.798 -1.5  2.798 1.5 /   
\plot -2.798 1.5  2.798 -1.5 /   
\plot -1.5 2.798  1.5 -2.798 /   
\endpicture
$$
$$
\beginpicture
\setcoordinatesystem units <1cm,1cm>         
\setplotarea x from -3 to 3, y from -3 to 3    
\put{$\omega_2$}[r] at -0.1 1.9
\put{$\alpha_2$} at -1.5 1.05
\put{$\alpha_1$}[l] at 1.05 0.15
\put{$\bullet$} at 0 0
\put{$\bullet$} at 1 0
\put{$\bullet$} at 2 0
\put{$\bullet$} at 3 0
\put{$\bullet$} at -1 0
\put{$\bullet$} at -2 0
\put{$\bullet$} at -3 0
\put{$\bullet$} at 0.5  0.8660
\put{$\bullet$} at 1.5  0.8660
\put{$\bullet$} at 2.5  0.8660
\put{$\bullet$} at -0.5  0.8660
\put{$\bullet$} at -1.5  0.8460
\put{$\bullet$} at -2.5  0.8660
\put{$\bullet$} at 0  1.732
\put{$\bullet$} at 1  1.732
\put{$\bullet$} at 2  1.732
\put{$\bullet$} at -1  1.732
\put{$\bullet$} at -2  1.732
\put{$\bullet$} at 0.5  -0.8660
\put{$\bullet$} at 1.5  -0.8660
\put{$\bullet$} at 2.5  -0.8660
\put{$\bullet$} at -0.5  -0.8660
\put{$\bullet$} at -1.5  -0.8660
\put{$\bullet$} at -2.5  -0.8660
\put{$\bullet$} at 0  -1.732
\put{$\bullet$} at 1  -1.732
\put{$\bullet$} at 2  -1.732
\put{$\bullet$} at -1  -1.732
\put{$\bullet$} at -2  -1.732
\linethickness=0.5pt                          
\plot 0 -0.04 0.72  -0.46 /    
\plot 0.72  -0.46 0.75  -0.43 /    %
\plot 0.75  -0.43 0 0  /    %
\plot 0.03 0 -0.44 -0.28  /  
\plot -0.44 -0.28 -0.23 -0.40 /  %
\plot -0.23 -0.40 -0.25 -0.44 /  %
\plot -0.25 -0.44 -0.52 -0.28 /  %
\plot -0.52 -0.28 -0.04 0 /  %
\setdots
\putrule from 0 3 to 0 -3        
\plot -1.5 -2.798   1.5 2.798 /  
\putrule from 3 0  to -3 0       
\plot -2.798 -1.5  2.798 1.5 /   
\plot -2.798 1.5  2.798 -1.5 /   
\plot -1.5 2.798  1.5 -2.798 /   
\endpicture
$$
These paths yield the following data:
$$\matrix{
\hbox{endpoint} &\hbox{maximal lift $\vec t$}  &\iota(\eta)=\tau_1
&v(w,\eta)^{-1} \cr
\omega_2 &(s_1) &\bar 1 &s_1 \cr
s_2\omega_2 &(s_2s_1) &\overline{s_2} &s_1s_2 \cr
s_2\omega_2-\alpha_1 &(s_1s_2s_1>s_2s_1) &\overline{s_1s_2} &s_1s_2 \cr
s_2\omega_2-2\alpha_1 &(s_1s_2s_1>s_2s_1) &\overline{s_1s_2} &s_1s_2 \cr
s_1s_2\omega_2 &(s_1s_2s_1) &\overline{s_1s_2} &s_1s_2s_1 \cr
s_2s_1s_2\omega_2 &(s_2s_1s_2) &\overline{s_2s_1s_2} &s_2s_1s_2 \cr
s_2s_1s_2\omega_2-\alpha_1 &(s_1s_2s_1s_2>s_2s_1s_2) &\overline{s_1s_2s_1s_2}
&s_2s_1s_2 \cr
s_2s_1s_2\omega_2-2\alpha_1 &(s_1s_2s_1s_2>s_2s_1s_2) &\overline{s_1s_2s_1s_2}
&s_2s_1s_2 \cr
s_1s_2s_1s_2\omega_2 &(s_1s_2s_1s_2) &\overline{s_1s_2s_1s_2} &s_2s_1s_2s_1 \cr
\alpha_1 &(s_2s_1s_2>s_1s_2>s_2) &\overline{s_2s_1s_2} &s_2 \cr
-\alpha_1 &(s_1s_2s_1s_2>s_2s_1s_2>s_1s_2) &\overline{s_1s_2s_1s_2} &s_2s_1 \cr
0 &(s_2s_1s_2>s_1s_2) &\overline{s_2s_1s_2} &s_2s_1 \cr
0 &(s_1s_2s_1s_2>s_2s_1s_2>s_1s_2>s_2) &\overline{s_1s_2s_1s_2} &s_2 \cr
}$$
and thus we get
$$\eqalign{
e^{\omega_2}T_{s_2s_1s_2s_1}
= 
T_{s_1}e^{\omega_2}
&+T_{s_1s_2}e^{s_2\omega_2}
+T_{s_1s_2}e^{s_2\omega_2-\alpha_1}
+T_{s_1s_2}e^{s_2\omega_2-2\alpha_1}
+T_{s_1s_2s_1}e^{s_1s_2\omega_2}  \cr
&
+T_{s_2s_1s_2}e^{s_1s_2s_1\omega_2}
+T_{s_2s_1s_2}e^{s_2s_1s_2\omega_2-\alpha_1}
+T_{s_2s_1s_2}e^{s_2s_1s_2\omega_2-2\alpha_1}  \cr
&+T_{s_2s_1s_2s_1}e^{s_1s_2s_1s_2\omega_2} 
+T_{s_2}e^{\alpha_1}
+T_{s_2s_1}e^{-\alpha_1}
+T_{s_2s_1}e^0
+T_{s_2}e^0
\cr
}$$
and
$$
e^{\omega_2}[\cO_{X_{s_1s_2s_1s_2}}]
=[\cO_{X_{s_1s_2s_1s_2}}]
+[\cO_{X_{s_1s_2s_1}}]
+3[\cO_{X_{s_2s_1s_2}}]
+3[\cO_{X_{s_2s_1}}]
+2[\cO_{X_{s_1s_2}}]
+2[\cO_{X_{s_2}}]
+[\cO_{X_{s_1}}].
$$

\section 5. Passage to $H^*(G/B)$

Let us explain how our results in $K(G/B)$ are related to the cohomology
$H^*(G/B)$ and Schubert polynomials.  The transfer is
by way of the Chern character $\ch$.

If $X$ is a finite CW complex then the {\it Chern character} (see [Mac, Ch. 10],
[Hi, \S23-24], [Ha, App. A])
$$\ch\colon K_{vb}(X)\otimes \QQ \longrightarrow H^*(X;\QQ)$$
is a natural ring isomorphism, i.e. if $f\colon X\to Y$ is continuous then
$\ch(f^!(x))=f^*(\ch(x))$.
If $f\colon X\to Y$ is a morphism of nonsingular projective varieties then
the Grothendieck-Riemann-Roch theorem [Ha, App. A Theorem 5.3], henceforth G-R-R,
says
$$\ch(f_!(x))=f_*(\ch(x)\td(\cT_f)),$$
where $\td(\cT_f)$ is Todd class of the relative tangent sheaf of $f$.
If $\cL$ is a line bundle on $X$ with first Chern class 
$\lambda\in H^2(X;\QQ)$ then the Chern character and the Todd class
of $\cL$ are the elements of $H^*(X)$ given by
$$\ch(\cL)=e^\lambda= \sum_{k\ge 0} {\lambda^k\over k!}\qquad\hbox{and}\qquad
\td(\cL)={\lambda\over 1-e^{-\lambda}},
\quad\hbox{respectively.}$$
The expression $e^\lambda$ is a finite sum since $\lambda^k=0$ in $H^*(X)$
whenever $k>\dim(X)$.

\medskip
\subsection{$H^*(G/B)$ as the quotient of a polynomial ring}

Let $X=G/B$.  Let $\gh$ be the Lie algebra of $T$ and
let $S(\gh^*)$ be the ring of polynomials on $\gh^*$ (over $\QQ$).  This is a
polynomial ring in the $n$ variables $\alpha_1,\ldots,\alpha_n$ (the simple
roots). Let $\hat\varepsilon\colon S(\gh^*)\to \QQ$ be the homomorphism given by
$\hat\varepsilon(\lambda)=0$ for all $\lambda\in \gh^*$.  If $f\in S(\gh^*)$ 
then $\hat\varepsilon(f)$ is the constant term of $f$.
It is a classical theorem of Borel (see [BGG, Prop. 1.3]) that 
$$H^*(G/B;\QQ) \cong {S(\gh^*)\over\hat\cI},\formula$$
where $\hat\cI$ is the ideal of $S(\gh^*)$ generated by
$\{f\in S(\gh^*)^W \ |\ f-\hat\varepsilon(f) = 0\}$. 
Proposition 1.5 is the K-theory analogue of Borel's theorem.

Let $-\lambda\in \Lambda$.
The element $-\lambda$ determines a character of $T$, denoted by
$e^{\lambda}\in R(T)$ (see section 1).  Let $c_1(\cL_{-\lambda})\in H^2(G/B)$
be the first Chern class of the line bundle $\cL_{-\lambda}=\phi(e^{\lambda})$
where $\phi\colon R(T)\to K(G/B)$ is the map in (1.2). 
Because of the isomorphism in (5.1) we often abuse notation and
write $\lambda=c_1(\cL_{-\lambda})\in H^*(G/B)$.  
All of the maps in the following
commutative diagram are isomorphisms. (Recall that we can identify $K(G/B)$ and
$K_{vb}(G/B)$.)
$$\matrix{
R(T)/\cI\otimes \QQ &\mapright{\phi} &K_{vb}(G/B)\otimes \QQ \cr
\mapdown{\ch} &&\mapdown{\ch} \cr
S(\gh^*)/\hat\cI &\mapright{\hat\phi} &H^*(G/B;\QQ) \cr
}
\qquad\hbox{where}\qquad
\matrix{
e^{\lambda} &\mapright{\phi} &[\cL_{-\lambda}] \cr
\mapdown{\ch} &&\mapdown{\ch} \cr
e^{\lambda} &\mapright{\hat\phi} &e^{c_1(\cL_{-\lambda})} \cr
}
\formula$$
The left hand $\ch$ map is obtained by viewing $R(T)$ and
$S(\gh^*)$ as subsets of $K(\cB_T)$ and $\widehat{H^*}(\cB_T)$, respectively,
where
$\cB_T$ is the classifying space of $T$.

\medskip
\subsection{The relation between $T_\alpha$ and $\partial_\alpha$}

Let $\alpha$ be a simple root and let 
$P_\alpha$ be the parabolic subgroup with Lie algebra $\gp_\alpha$
spanned by $\gb$ and the root space $\gg_{-\alpha}$.  Let $f_\alpha\colon G/B\to
G/P_\alpha$ corresponding $\PP_1$-bundle.

\medskip\noindent
{\it A proof of Proposition 3.3:}\ \ 
Since $K(G/B)$ is torsion free, we can check
the identity $f_\alpha)^!(f_\alpha)_!(x)=T_\alpha(x)$ by applying $\ch$ to both
sides and checking the result in cohomology.

The G-R-R for the map $f_\alpha$ says
$$\ch((f_\alpha)_!(x))=(f_\alpha)_*(\ch(x)\td(\cT_{f_\alpha})),
\formula$$
where $\cT_{f_\alpha}$ is the bundle of tangents along the fibres, which is
the line bundle associated to the $\ad(B)$-module $\gp_\alpha/\gb$.  Since this
module has weight $-\alpha$, $c_1(\cT_{f_\alpha})=\alpha\in H^2(G/B\colon \ZZ)$
and so (5.3) becomes
$$\ch((f_\alpha)_!(x))=(f_\alpha)_*\left(\ch(x){\alpha\over 1-e^{-\alpha}}\right).
\formula$$
The BGG operator $\partial_\alpha=(f_\alpha)^*(f_\alpha)_*$ is explicitly given
by the formula
$$\partial_\alpha(z) = {z-s_\alpha(z)\over \alpha},
\formula
$$
see [D].  Thus, by applying $(f_\alpha)^*$ to
both sides of (5.4) we obtain
$$\ch((f_\alpha)^!(f_\alpha)_!(x))
=(f_\alpha)^*(f_\alpha)_*\left(\ch(x){\alpha\over 1-e^{-\alpha}}\right)
=\partial_\alpha\left(\ch(x){\alpha\over1-e^{-\alpha}}\right).
\formula
$$
The strategy now is to manipulate the right hand side of (5.6) using the
``skew-Leibniz'' rule satisfied by $\partial_\alpha$ to obtain
$\ch(T_\alpha(x))$.  For convenience, let $y=\ch(x)$.  Then the right side of
(5.6) is
$$\partial_\alpha\left(y{\alpha\over 1-e^{-\alpha}}\right)
=y\partial_\alpha \left({\alpha\over 1-e^{-\alpha}}\right)+
s_\alpha\left({\alpha\over 1-e^{-\alpha}}\right)\partial_\alpha(y)
$$
and we claim
$$
\hbox{(a)}\quad s_\alpha\left({\alpha\over
1-e^{-\alpha}}\right)= {\alpha\over e^\alpha-1},
\qquad\qquad
\hbox{(b)}\quad \partial_\alpha \left({\alpha\over 1-e^{-\alpha}}\right)
=1.
$$
The first equality is trivial and the second can be proved by formal computation
(carefully done!) or by applying the G-R-R again.  
Using (a) and (b) we obtain
$$\partial_\alpha\left(y{\alpha\over 1-e^{-\alpha}}\right)
=y+{\alpha\over e^\alpha-1}\left({y-s_\alpha(y)\over\alpha}\right).$$
Now cancelling the $\alpha$'s in the second term on the right and recalling that
$y=\ch(x)$ we find
$$\ch((f_\alpha)^!(f_\alpha)_!(x))
=\ch(T_\alpha(x)). \qquad\hbox{\qed}
$$

We see that 
$$\ch(T_\alpha(x))
=\partial_\alpha\left(\ch(x) {\alpha\over 1-e^{-\alpha}}\right)
\formula$$
which relates $T_\alpha$ to $\partial_\alpha$ in an explicit way.  In fact, if
one wished one could reverse the argument and derive the formula (5.5)
for $\partial_\alpha$ from the formula for $T_\alpha$.

\subsection{A proof of Proposition 1.6}

\bigskip\noindent{\bf Proposition 1.6. }\sl
If $f\colon G/B\to G/P$ is the natural projection then the induced map
$f^!\colon K(G/P)\to K(G/B)$ is an injection.
\pf  Since $K(G/P)$ is torsion free there is a natural injection
$K(G/P)\hookrightarrow K(G/P)\otimes \QQ\to H^*(G/P;\QQ)$ and so it suffices to
check that the pull-back $f^*$ in rational cohomology is injective.  Since the
odd cohomology groups of the base and fiber are zero the Serre
spectral sequence of the bundle $P/B\to G/B\to G/P$ shows that 
$f^*$ is injective even for $H^*(G/P;\ZZ)$.
\endpf

\vfill\eject

\subsection{Dictionary between $K(G/B)$ and $H^*(G/B)$}

In summary, the Chern character gives an isomorphism
$$\matrix{
R(T)/\cI &\cong &K(G/B)  &\mapright{\ch}  &H^*(G/B) &\cong &S(\gh^*)/\hat\cI\cr
&&e^\lambda &\longmapsto &e^\lambda,\cr
}$$
where
$$\matrix{
\cI=\hbox{ideal generated by $\{f\in R(T)^W\ |\ f-\varepsilon(f)=0\}$},
&\qquad 
&\matrix{
\varepsilon\colon &R(T)&\longrightarrow &\ZZ \cr
&e^\lambda &\longmapsto &1,\cr
} \cr
\cr\cr
\hat\cI 
=\hbox{ideal generated by $\{f\in S(\gh^*)^W\ |\ f-\hat\varepsilon(f)=0\}$},
&&\matrix{
\hat\varepsilon\colon &S(\gh^*)&\longrightarrow &\ZZ \cr
&\lambda &\longmapsto &0.\cr
} \cr
}$$
Let $[X_w]\in H^*(G/B)$ be the element which is Poincar\'e dual
to the  fundamental cycle of $X_w$ in $H_*(G/B)$.  This
element is called a {\it Schubert polynomial}.
Then 
$$\hbox{
$K(G/B)$ has basis
$\{[\cO_{X_w}]\ |\ w\in W\}$}\qquad\hbox{and}\qquad 
\hbox{$H^*(G/B)$ has basis $\{ [X_w]\ |\ w\in W\}$.}
$$
From a general fact  (see [Fu, Ex. 15.2.16] or [CG, 5.8.13(i) and p. 289])
$$\ch([\cO_{X_w}])=[X_w]+\hbox{higher degree terms}.$$
where $\deg([X_w])=\dim(G/B)-\dim(X_w)=N-\ell(w)$, where $N$ is the
number of positive roots for $\gg$.

If $\alpha$ is a simple root and $f_\alpha\colon G/B\to G/P_\alpha$ is the
corresponding $\PP^1$-bundle, then
$$T_\alpha(x) = (f_\alpha)^!(f_\alpha)_!(x)
={e^\alpha x -s_\alpha(x)\over e^\alpha-1}
\qquad\hbox{and}\qquad
\partial_\alpha(x) =(f_\alpha)^*(f_\alpha)_*(x)
={x-s_\alpha(x)\over \alpha}$$
in $K(G/B)$ and $H^*(G/B)$, respectively.  As illustrated in (5.7)
each of these two formulas can be derived from the other via the use of the
Grothendieck-Riemann-Roch Theorem.  This means that the following
formulas (Proposition 3.4 and [BGG, Th. 3.14])
$$T_\alpha([\cO_{X_w}])=\cases{
[\cO_{X_{ws_\alpha}}], &if $ws_\alpha>w$,\cr
[\cO_{X_w}], &if $ws_\alpha<w$,}
\qquad\hbox{and}\qquad
\partial_\alpha([X_w])=\cases{
[X_{ws_\alpha}], &if $ws_\alpha>w$,\cr
[X_w], &if $ws_\alpha<w$,} 
$$
are equivalent.

Our new Pieri-Chevalley formula in $K(G/B)$, Theorem 4.3, and Chevalley's
classical Pieri formula in $H^*(G/B)$, [Ch, Prop. 10], are
$$e^\lambda [\cO_{X_w}] = 
\sum_{\eta\in \cT^\lambda_w} [\cO_{X_{v(\eta,w)}}] 
\qquad\hbox{and}\qquad
\lambda\cdot [X_w] = \sum_{{}\atop v\mapright{\alpha}w}
\langle\lambda,\alpha^\vee\rangle [X_v],$$
where the sum is over all $v\in W$ such that $\ell(v)=\ell(w)-1$ and there
is a root $\alpha$ such that $v=s_\alpha w$.
Chevalley's formula can be obtained from ours formula by subtracting
$[\cO_{X_w}]$ from each side, applying the Chern character $\ch$, and comparing
the lowest degree terms on each side.

\vfill\eject
\bigskip\bigskip
\centerline{\smallcaps References}



\medskip
\item{[BGG]} {\smallcaps I.N. Bernstein, I.M. Gel'fand and S.I. Gel'fand},
{\it Schubert cell and cohomology of the spaces $G/P$}, Russ. Math. Surv. {\bf 28}
(3) (1973), 1--26.


\medskip
\item{[BS]} {\smallcaps A. \ Borel and J.-P.\ Serre},
{\it Le Th\'eor\`eme de Riemann-Roch (d'apres Grothendieck)},
Bull. Soc. Math. France {\bf 68} (1958), 97--136.

\medskip
\item{[C]} {\smallcaps C.\ Chevalley}, {\it Le classes d'equivalence
rationelle}, I, II, Seminaire C. Chevalley, Anneaux de Chow et applications
(mimeographed notes), Paris, 1958.

\medskip
\item{[Ch]} {\smallcaps C.\ Chevalley}, {\it Sur les decompositions cellulaires
des espaces $G/B$},  in {\sl Algebraic Groups and their Generalizations:
Classical Methods}, W. Haboush and B. Parshall eds.,
Proc. Symp. Pure Math., Vol. {\bf 56} Pt. 1, Amer. Math. Soc. (1994), 1--23.

\medskip
\item{[CG]} {\smallcaps N.\ Chriss and V.\ Ginzburg}, 
{\sl Representation theory and complex geometry}, Birkh\"auser, Boston, 1997.

\medskip
\item{[D]} {\smallcaps M.\ Demazure}, 
{\it D\'esingularisation des vari\'et\'es de
Schubert g\'en\'eralis\'ees}, Ann. Sci. \'Ecole Norm. Sup. {\bf 7} (1974), 53--88.

\medskip
\item{[FL]} {\smallcaps W.\ Fulton and A.\ Lascoux},
{\it A Pieri formula in the Grothendieck ring of a flag bundle},
Duke Math. J. {\bf 76} (1994), 711--729.

\medskip
\item{[Fu]} {\smallcaps W.\ Fulton},
{\sl Intersection Theory},  Ergebnisse der Mathematik (3) {\bf 2},
Springer-Verlag, Berlin-New York, 1984. 

\medskip
\item{[Ha]} {\smallcaps R.\ Hartshorne}, {\sl Algebraic geometry}, 
Graduate Texts in Mathematics {\bf 52}, Springer-Verlag, New York-Heidelberg,
1977.

\medskip
\item{[H]} {\smallcaps F.\ Hirzebruch}, {\sl Topological methods in algebraic
geometry}, Classics in Mathematics, Springer-Verlag, Berlin, 1995. 

\medskip
\item{[Ke]} {\smallcaps G.\ Kempf}, {\it Linear systems on homogeneous spaces},
Ann. Math. {\bf 103} (1976), 557--591.

\medskip
\item{[KK]} {\smallcaps B.\ Kostant and S.\ Kumar}, 
{\it $T$-equivariant K-theory of generalized flag varieties}, J. Differential
Geom. {\bf 32} (1990), 549--603. 

\medskip
\item{[L1]} {\smallcaps P.\ Littelmann},
{\it A Littlewood-Richardson rule for symmetrizable Kac-Moody algebras},
Invent. Math. {\bf 116} (1994), 329-346.

\medskip
\item{[L2]} {\smallcaps P.\ Littelmann},
{\it Paths and root operators in representation theory}, Ann. Math.
{\bf 142} (1995), 499-525.

\medskip
\item{[L3]} {\smallcaps P.\ Littelmann},
{\it Characters of representations and paths in ${\goth H}_{\RR}^*\;$},
Proc. Symp. Pure Math. {\bf 61} (1997), 29-49.

\medskip
\item{[Mac]} {\smallcaps I.G.\ Macdonald},
{\sl Algebraic geometry: Introduction to schemes}, W. A. Benjamin, New
York-Amsterdam,1968.

\medskip
\item{[P]} {\smallcaps H.\ Pittie}, 
{\it Homogeneous vector bundles over homogeneous spaces},
Topology {\bf 11} (1972), 199--203.

\medskip
\item{[Ra]} {\smallcaps A. Ramanathan}, 
{\it Frobenius splitting and Schubert varieties}, in {\sl Proceedings of the
Hyderabad Conference on Algebraic Groups} (Hyderabad, 1989),  Manoj
Prakashan, Madras, 1991, p. 497--508. 

\medskip
\item{[W]} {\smallcaps A.\ Weil}, 
{\it Demonstration topologique d'un th\'eor\`eme fondamental de Cartan},
C.R. Acad. Sci. {\bf 200} (1935), 518--520.

\vfill\eject
\end

 The
BGG-operator is the operator
$\partial_\alpha$ on $H^*(G/B)$ given by
$$\partial_\alpha(x)=(f_\alpha)^*(f_\alpha)_*(x)$$
(see [D]).
The operators $T_\alpha$ from section 3 are the K-theoretic
analogues of the $\partial_{\alpha}$.  We can use the formula defining
the $T_\alpha$ operators and the Grothendieck-Riemann-Roch theorem to derive
a formula for the $\partial_\alpha$ operators as follows.

Let $X=G/B$ and $Y=G/P_\alpha$. Since $f_\alpha\colon X\to Y$
is a fiber bundle,
$$\td(\cF_X)=(f_\alpha)^*(\td(\cF_Y))\td(\cF),$$
where $\cF$ is the bundle of tangents along the fibers.
Since $(f_\alpha)_*(\td(\cF_Y))=\td(\cF_Y)$ the
Grothendieck-Riemann-Roch theorem gives
$$\ch((f_\alpha)_!(x))\td(\cF_Y) =
(f_\alpha)_*(\ch(x)\td(\cF))\td(\cF_Y).$$ 
Since $\td(\cF_Y)$ is a power series with constant term $1$,
it is an invertible element of cohomology and we can cancel it from both
sides.  Now apply $f_\alpha^*$ to both sides
to get
$$\ch((f_\alpha)^!(f_\alpha)_!(x)) = (f_\alpha)^*(f_\alpha)_*(\ch(x)\td(\cF)).$$
Since the bundle $\cF$ of tangents along the fibers of 
$f_\alpha\colon G/B\to G/P_\alpha$ is the line bundle $\cL_{-\alpha}$
(see [Ke, \S 1 Lemma 1]) we have
$$T_\alpha(x) = \ch((f_\alpha)^!(f_\alpha)_!(x)) 
= (f_\alpha)^*(f_\alpha)_*\left(\ch(x){\alpha\over 1-e^{-\alpha}}\right) 
=\partial_\alpha\left(x{\alpha\over 1-e^{-\alpha}}\right).$$
(By the diagram in (5.2) we can abuse notation in cohomology and write $x$
and $T_\alpha(x)$ for
$\ch(x)$ and  $\ch((f_\alpha)^!(f_\alpha)_!(x))$, respectively.)
This identity allows us to use the definition of $T_\alpha$ from section
3 and derive the classical formula [BGG] for $\partial_\alpha$ from
the formula for $T_\alpha$ given in section 3:
$$
\partial_\alpha(x)
=\partial_\alpha
\left(x{1-e^{-\alpha}\over\alpha}{\alpha\over 1-e^{-\alpha}}\right)
=T_\alpha\left(x{1-e^{-\alpha}\over\alpha}\right)
={x-s_\alpha(x)\over \alpha}.\formula
$$
The factor $(1-e^{-\alpha})/\alpha$ makes sense since the Todd class
$\alpha/(1-e^{-\alpha})$ of the line bundle $\cL_{-\alpha}$ is an invertible
element of cohomology.  Of course one could reverse this argument
to use the classical [BGG] formula for $\partial_\alpha$ to derive the
formula for $T_\alpha$ given in section 3.

\medskip
\subsection{Recovering Chevalley's formula from Theorem ???}

Given any identity in $K(G/B)$ we can apply $\ch$ to both sides to
obtain an identity in $H^*(G/B)$.  Since $H^*(G/B)$ is a graded ring we
thus obtain a cohomology identity in every degree.  Let $[X_w]$ be the element of $H^*(G/B)$ which is Poincar\'e dual to the 
fundamental cycle of $X_w$ in $H_*(G/B)$.
Then 
$$\hbox{
$K(G/B)$ has basis
$\{[\cO_{X_w}]\ |\ w\in W\}$}\qquad\hbox{and}\qquad 
\hbox{$H^*(G/B)$ has basis $\{ [X_w]\ |\ w\in W\}$.}
$$
From a general fact  (see [Fu, Ex. 15.2.16] or [CG, 5.8.13(i) and p. 289])
$$\ch([\cO_{X_w}])=[X_w]+\hbox{higher degree terms}.$$
Furthermore, $\deg([X_w])=\dim(X_w)=N-\ell(w)$, where $N$ is the dimension
of $G/B$.

We can rewrite Theorem ??? as
$$(e^\lambda-1)[\cO_{X_w}]
= \sum_{
\iota(\eta)\ne \bar w} [\cO_{X_{v(\eta,w)}}],$$
where $\bar w=wW_\lambda\in W/W_\lambda$ and the sum is over all paths $\eta\in
\cT^\lambda_w$ such that $\iota(\pi)\ne \bar w$.
If we apply $\ch$ to both sides and compare lowest degree terms 
we obtain
$$\lambda\cdot [X_w]
= \sum_{\alpha\atop v(\eta,w)\longrightarrow w\ \ \ \ } [X_{v(\eta,w)}],$$
where the sum is over all paths $\eta\in \cT^\lambda_w$ such that there
is a positive root $\alpha$ with $s_\alpha v(\eta,w) = w$ and
$\ell(w)=\ell(v(\eta,w))+1$.
This last formula is equivalent to Chevalley's formula
in $H^*(G/B)$, see ??? below.

\prop If $f\colon G/B\to G/P$ is the natural projection then the induced map
$f^!\colon K(G/P)\to K(G/B)$ is an injection.
\pf  For any finite CW
complex $X$, there is an isomorphism $\ch\colon K^*(X)\otimes \QQ\to
H^*(X;\QQ)$.  Since $K(G/P)$ is torsion free there is a natural injection
$K(G/P)\hookrightarrow K(G/P)\otimes \QQ\to H^*(G/P)$ and so it suffices to
check that the pull-back $f^*$ in rational cohomology is injective.  In fact,
$f^*$ is injective even for $H^*(\ ;\ZZ)$ as one sees easily from the Serre
spectral sequence of the bundle $P/B\to G/B\to G/P$; because the odd cohomology
groups of the base and fiber are zero.
\endpf

We can use the Grothendieck-Riemann-Roch theorem to give alternative proofs
of several of our earlier results.

\prop If $f\colon G/B\to G/P$ is the natural projection then the induced map
$f^!$ is an injection.
\pf Since $f\colon G/B\to G/P$ is the projection of an
algebraic fiber bundle, the Grothendieck Riemann-Roch theorem gives
$$\ch(f_!(x))=f_*(\ch(x)\td(\cF)),$$
where $\cF$ is the bundle of tangents along the fibers.  Taking $x=1$,
$$f_*(\td(\cF))=\td(T(P/B))[P/B]=1,$$
by the Hirzebruch-Riemann-Roch theorem, because $P/B$ is a rational
variety.  Since $\ch$ is injective, we have $f_!(1)=1$ and so
$f_!(f^!(x))=x$ for every $x\in K(G/P)$.  That is, $\Im f^!$ is an 
additive summand of $K(G/B)$. 
\endpf

\prop  Let $\alpha_j$ be a simple root and let $T_j=T_{\alpha_j}$.
For every $x\in K(G/B)$, 
$$(f_j)^!\circ(f_j)_!(x)=T_j(x).$$
\endprop
\pf  Let us first see that the bundle in ??? is the projectivization
of a rank 2 vector bundle over $G/P_j$. 
Let $\omega_j$ be the $j$th fundamental weight and let $V^{\omega_j}$ be the
irreducible $G$-module with highest weight 
$\omega_j$ and let $E_j$ be the two-dimensional $P_j$-submodule of
$V^{\omega_j}$ spanned by the highest weight vector $v^+$ and $X_{-\alpha_i}v^+$
where $X_{-\alpha_i}$ is a nonzero element of $\gg_{-\alpha_j}$.  The weights of
$E_j$ are $\omega_j$ and $\omega_j-\alpha_j$. Consider the vector
bundle ${\cal E}=G\times_{P_j} E_j\to G/P_j$.  One can easily construct a bundle
isomorphism
$$G\times_{P_j} \PP_1(E_j)\mapright{\cong} G/B$$
and under this identification we obtain a line sub-bundle ${\cal
L}_{\omega_j}\subseteq f^*_j({\cal E}_j)$, such that ${\cal L}_{\omega_j}$
restricted to any fiber $\PP_1$ is $\cO_{\PP_1}(-1)$.

Since $G/B\to G/P_j$ is the $\PP_1$-bundle of $\cE_j$, Bott periodicity implies 
(see [A], section ???) that $K(G/B)$ is a free $K(G/P_j)$-module on two
generators, $1$ and $e^{-\omega_j}$.  The map 
$T_j$ is a $K(G/P_j)$-module map
since $K(G/P_j)=K(G/B)^{W(P_j)}$,
and $(f_j)^!\circ (f_j)_!$ is a $K(G/P_j)$-module map
because for any map of varieties $f\colon X\to Y$,
$f_!(f^!(y)) =yf_!(1)$.  This means that it suffices to check that
$$T_j(1)=(f_j)^!\circ(f_j)_!(1)\qquad\hbox{and}\qquad
T_j(e^{-\omega_j})=(f_j)^!\circ(f_j)_!(e^{-\omega_j}).$$  
The first equality follows from $(f_j)_!(1)=1$ and $T_j(1)=1$.

Choose an open  set $U\subseteq G/P_j$ over which the bundle $G/B \to G/P_j$ is
trivial so that $f_j^{-1}(U)\cong U\times \PP_1$.  Then
$$R^q(f_j)_!(\cL_{\omega_j})(U)=H^q(U\times \PP_1,\cL_{\omega_j}).$$
Since $\cL_{\omega_j}$ is constant in the $U$-directions, one can
use Grothedieck's Kunneth formula [BS] to obtain
$$H^q(U\times \PP_1,\cL_{\omega_j})
= \sum H^i(U,\cO_U)\otimes H^j(\PP_1,\cO_{\PP_1}(-1)).$$
Since $H^j(\PP_1,\cO(-1))=0$ for all $j\ge 0$ we get that
$$R^q(f_j)_!(\cL_{\omega_j})=0,
\qquad\hbox{for all $q$, including zero.}$$
Thus $(f_j)^!\circ(f_j)_!(e^{-\omega_j})=0$.  
Finally, a direct computation shows that
$T_j(e^{-\omega_j})=0$. 
\endpf